\documentclass[11pt,a4paper]{article}
\usepackage{amsfonts}
\usepackage{mathrsfs}

\newtheorem{prop}{Proposition}
\newtheorem{lemma}[prop]{Lemma}
\newtheorem{thm}[prop]{Theorem}
\newenvironment{demo}{\par\noindent {\em Proof}.\/}{
\hfill $\Box$ \medskip\\}
\newenvironment{rmk}{\addtocounter{prop}{1}%
\medskip\par\noindent\textbf{Remark \theprop\ }}{\medskip\par}

\renewcommand{\P}{\mathbb{P}}
\newcommand{\E}{\mathbb{E}}

\newcommand{\erre}{\mathbb{R}}

\renewcommand{\epsilon}{\varepsilon}

\newcommand{\ds}{\displaystyle}

\title{The Stochastic Goodwill Problem}

\author{Carlo Marinelli\\[4pt]
\small\em Institut f\"ur Angewandte Mathematik, Universit\"at Bonn\\
\small\em Wegelerstr. 6, D-53115 Bonn (Germany)\\
\small{\em e-mail:} \texttt{cm1310@gmail.com}\\[10pt]}

\date{\normalsize August 10, 2005}

\begin{document}
\maketitle

\begin{abstract}\noindent
  Stochastic control problems related to optimal advertising under
  uncertainty are considered. In particular, we determine the optimal
  strategies for the problem of maximizing the utility of goodwill at
  launch time and minimizing the disutility of a stream of advertising
  costs that extends until the launch time for some classes of
  stochastic perturbations of the classical Nerlove-Arrow dynamics. We
  also consider some generalizations such as problems with constrained
  budget and with discretionary launching.
\bigskip\par
\noindent{\it Key Words:} optimal advertising under
uncertainty, Bellman equation, control with discretionary stopping.
\end{abstract}

\section{Introduction}
We consider the optimization problem faced by a firm that, while
advertising a product prior to its introduction to the market, wants
to determine the optimal advertising policy for the maximization of
the product image (also called \emph{goodwill}), and the minimization
of the total discounted cost.  We shall also consider the problem of
optimizing the launching time, thus allowing the firm to decide at its
discretion to stop the advertising campaign and start selling the
product.

This type of problems can be traced back at least to
Nerlove and Arrow \cite{NA}, who proposed to model the stock of advertising
goodwill $x(t)$ at time $t\geq 0$ as
\begin{equation}
\label{eq:NA}
\dot{x}_t = u_t - \rho x_t, \quad x_0=x\geq 0,
\end{equation}
where $u_t$ is the rate of advertising expenditure, $\rho>0$ is a
factor of deterioration of product image in absence of advertisement.
The optimization problem for a firm that seeks to maximize awareness
of its product at a given launch time $T>0$ and to minimize its
advertising effort until $T$ could be formulated as a
multi-objective program of the type
\begin{equation}
\label{eq:mob}
\left(
\sup_{u\in\mathcal{U}} \E[e^{-cT}\varphi(x_T)],
\inf_{u\in\mathcal{U}} \E\Big[\int_0^T e^{-ct}h(u_t)\,dt\Big]
\right),
\end{equation}
subject to the dynamics (\ref{eq:NA}), where $c>0$ is a discount
factor, $\varphi:\erre\to\erre$ is a reward function,
$h:\erre\to\erre$ is a loss function,
and $\mathcal{U}$ is the set of measurable functions $u:[0,T]\to U$,
with $U$ a closed subset of $\erre_+$. Following a standard procedure
in multi-objective optimization (see e.g. Zeleny \cite{zeleny}), one
takes a weighted average of the objectives in (\ref{eq:mob}) and
obtains the deterministic optimal control problem
\begin{equation}
\label{eq:det}
\sup_{u\in\mathcal{U}}
\left( e^{-c T}\gamma_0\varphi(x_T) -
\int_0^T e^{-c t}h(u_t)\,dt\right),
\end{equation}
where $\gamma_0$ is a positive constant.
Sufficient conditions for the problem to be well posed are, e.g., that
$\varphi$ is concave and continuous, $h$ is convex and continuous, and
$U$ is compact. The conditions are also meaningful from an economic
point of view, as it is customary to use concave (increasing) utility
functions as measures of reward, and convex (increasing) loss functions.
In the simplest case,
one could take $h(u)=u$, so that the second term in (\ref{eq:det})
coincides with the total discounted advertising expenditure.
This deterministic optimal control formulation has been extended by
many authors to account for delay effects, non-linearity in the
response to advertisement, and many other factors. For a recent
work on the subject, which also contains a list of related references,
we refer to Buratto and Viscolani \cite{BV}.

On the other hand, less work has been devoted to the case of
stochastic evolution of goodwill level: for a few examples of works in
this direction, we refer to the survey by Feichtinger, Hartl and Sethi
\cite{sethi} and references therein, and to the more recent papers of
Grosset and Viscolani \cite{GV} and Buratto and Grosset \cite{BG}.
The emergence of randomness in the dynamics of goodwill is quite
natural for several reasons: one may think, for example, that random
fluctuations in the goodwill level are the effect of external factors
beyond the control of the firm, or that noise enters through the
control, since the effect of advertising may be partly uncertain (see
section \ref{sec:models} for a detailed discussion).

In this work we introduce some stochastic extensions of the classical
model of Nerlove and Arrow and study related optimization problems.
We do not aim at maximum generality, instead we focus on models whose
special structure allows us to obtain explicit solutions.  In
particular, in section \ref{sec:models} we propose a stochastic
extension of the Nerlove-Arrow dynamics and motivate it by marketing
assumptions. We formulate a rather general problem of optimizing an
objective function that weighs (a function of) product image at a
fixed time and the cumulative cost of advertising effort, and we
construct a nearly optimal advertising strategy.
The special case of linear reward of goodwill and linear cost of advertising
effort is considered in section 3: the special structure of the problem
allows one to obtain the value function and the optimal policy in closed
form, and to consider more general problems of advertising with a limited
budget.
In section 4 we study another case where explicit solutions can be obtained,
i.e. the case of quadratic reward of final goodwill and quadratic cost of
advertising.
Under these assumptions we also explicitly solve in section 5 a problem of
optimal advertising with discretionary stopping to reach a target level of
product awareness.
We conclude suggesting some problems not addressed in this paper.

\section{Stochastic models for goodwill dynamics and related optimization problems}
\label{sec:models}
Let $x_t$ be the level of product image at time $t$, $0 \leq t\leq T$,
where $T>0$ is the end of the planning horizon (the time at which the
product will be launched).  We postulate a dynamics for $x_t$ given by
the following stochastic differential equation, which is a stochastic
perturbation of the Nerlove-Arrow dynamics (\ref{eq:NA}):
\begin{equation}
\label{eq:dyn}
dx_t = (- \rho x_t + u_t)\,dt + \sigma(x_t,u_t)\,dw_t, \quad x_0=x\geq 0,
\end{equation}
where $\rho$ is a positive constant, $\sigma:\erre^2\to\erre$ is
Lipschitz continuous, and $w_t$ is a standard real valued Brownian
motion on a filtered probability space
$(\Omega,\mathcal{F},\mathbb{F},\P)$,
$\mathbb{F}=(\mathcal{F}_t)_{t\in[0,T]}$.
The control process $u_t$ models the rate of advertising \emph{effort} by
the firm, and is assumed to be measurable, adapted, and taking values in a
closed convex subset $U$ of $[0,+\infty)$. We will denote by $\mathfrak{U}$
the set of controls satisfying these properties. We use gross rating points
(GRPs) to measure advertising effort, instead of the rate of advertising
expenditure, following a recent trend in the marketing literature -- see
e.g. Dube and Manchanda \cite{dube}, Vilcassim, Kadiyali and Chintagunta
\cite{vilca}.

As briefly mentioned in the introduction, there are several reasons to
study stochastic extensions of the Nerlove-Arrow model. Early papers
on the subject such as Rao \cite{rao} and Raman \cite{raman} advocated
the use of stochastic models with the observation that such effects as
copy and competitive changes would render uncertain the effect of
advertising on goodwill. In both papers stochastic perturbations of
the Nerlove-Arrow dynamics of the type (\ref{eq:dyn}) were proposed,
with $\sigma$ not depending on $x_t$ nor $u_t$.

While a model with only additive noise can be useful as a first
approximation, it is reasonable to consider more general stochastic
disturbances, depending also on the goodwill level $x_t$ and on the
intensity of advertising effort $u_t$. In particular, one could
distinguish, potentially among others, three sources of uncertainty: a
``background noise'', of additive type, due to the kind of effects
indicated by Raman \cite{raman} and Rao \cite{rao}, which are not
directly influenced neither by the popularity of the product nor by
the intensity of advertising. A second contribution to the intensity
of noise can be attributed to the uncertainty in the opinion of the
(potential) customers that are aware of the product.  Finally, one
should take into account the uncertainty in the effect of advertising.
In order to model explicitly such disturbances, it may be reasonable
to assume that the intensities of the second and third type of noise
just mentioned are proportional respectively to the goodwill level
$x_t$ (a proxy for the number of customers that are aware of the
advertised product) and to the level of advertising effort $u_t$.
In particular, the noise component proportional to the goodwill level could
also be interpreted as the effect of ``internal influence'' (also called
word-of-mouth communication), due to the random outcome on the goodwill
level of the interaction between customers who know the product. Similar
ideas (and terminology) are extensively employed in the literature on new
product diffusion (see e.g. Bass \cite{bass}).
The third source of uncertainty described above could be suggestively
justified by attaching to each unit of advertising effort a random effect,
so that, heuristically speaking, ``noise enters the system through the
control'' via $u_t \mapsto u_t(1+\sigma_2\frac{dw_t}{dt})$.

We would also like to mention that models for the evolution of
goodwill expressed as stochastic perturbations of the Nerlove-Arrow
dynamics of the type (\ref{eq:dyn}), where the diffusion coefficient
$\sigma$ depends both on $x_t$ and on $u_t$, appeared also as
diffusion approximations of models based on discrete-time Markov
chains (see e.g. Tapiero \cite{tapiero} and references therein), and
in the above mentioned papers \cite{GV}, \cite{BG}.

Due to the lack of empirical studies and of theoretical papers in the
marketing literature on the determinants of uncertainty in the
dynamics of goodwill, we are led to consider, in the same spirit of
the cited works of Raman and Rao, a linear specification of the
diffusion coefficient $\sigma$:
$$
\sigma(x,u) = \sigma_0 + \sigma_1 |x| + \sigma_2 u,
$$
where $\sigma_0$, $\sigma_1$, $\sigma_2$ are fixed non-negative
constants. In practice, the values of these coefficients should be
determined by ad hoc empirical studies and/or by specific managerial
and marketing insights.

A natural analog in the stochastic setting of the general optimization
problem (\ref{eq:det}) can now be formulated. Let us define the
performance functional relative to strategy $u \in \mathfrak{U}$ as
\begin{equation}
\label{eq:fo}
v^u(s,x) = \E^u_{s,x}\left[
e^{-cT}\varphi(x_T) - \int_0^T e^{-ct} h(u_t)\,dt
\right],
\end{equation}
and the value function as
\begin{equation}
\label{eq:value}
v(s,x)=\sup_{u\in \mathfrak{U}} v^u(s,x),
\end{equation}
where $h:U\to\erre_+$ is bounded and $|\varphi(x)|<K(1+|x|^m)$ for
some positive constants $K$, $m$. By $\E^u_{s,x}$ we mean, as usual,
expectation with respect to the law of the controlled diffusion
$$
x_t = x + \int_s^t (-\rho x_r + u_r)\,dr + 
\int_s^t (\sigma_0 + \sigma_1 |x_r| + \sigma_2 u_r)\,dw_r,
\quad s \leq r \leq t \leq T.
$$
Our objective is to characterize the value function and to find (or
approximate) strategies realizing the supremum in (\ref{eq:value}).
Although the problem is well posed under the given assumptions,
particularly meaningful choices from the economic point of view are
$\varphi$ concave increasing (we have in mind the utility function of a
risk-averse agent) and $h$ convex increasing (typical choice of cost
function).

We shall study the problem through the dynamic programming approach,
i.e. through the study of the associated Bellman equation, which can be
written as
\begin{equation}
\label{eq:bellman}
\ds \sup_{u\in U}
\Big[\frac{\partial\psi}{\partial t} + L^u\psi - c\psi - h(u) \Big] = 0, \quad
\psi(T,x) = \varphi(x),
\end{equation}
where $L^u$ is the differential operator defined by
$$
L^u = \frac12a(x,u)\partial_x^2 + b(x,u)\partial_x,
$$
with $a:=\sigma\sigma^*=\sigma^2$ and $b(x,u)=-\rho x+u$.

Equation (\ref{eq:bellman}) is (under the assumption $\sigma_0\neq 0$)
a fully nonlinear uniformly nondegenerate parabolic PDE, for which
general results about existence of smooth solutions are available
under additional assumptions of smoothness and boundedness of the
coefficients and of the reward and loss functions. However, we can
prove existence of nearly optimal control in a general setting that
covers also the case $\sigma_0=0$, for which the Bellman equation
(\ref{eq:bellman}) becomes degenerate. In particular we have the
following result.
\begin{thm}
  For any $\epsilon>0$, $s\in[0,T]$, $x\in\erre$, there exists a
  Markov strategy $u^\epsilon_t \in \mathfrak{U}$ such that $v(s,x)
  \leq v^{u^\epsilon}(s,x)+\frac43\epsilon$.
\end{thm}
\begin{demo}
We divide the proof in three steps. In the first step we introduce a
sequence of approximating problems, in the second step we construct
an optimal feedback control for each approximating problem. In the
last step we prove that optimal controls for the approximating
problems are nearly optimal for the original problem.
\par\noindent
\textsc{Step 1:}
Let $\zeta_n \in C^\infty_0(\erre)$ with $\int_\erre
\zeta_n(x)\,dx=1$, $n=0,1,2\ldots$, be a sequence of mollifiers.
Define
$$
\tilde{b}_n(x,u) = \left\{
\begin{array}{ll}
\ds \rho n + u, & x<-n\\
\ds -\rho x + u, & -n \leq x \leq n\\
\ds -\rho n + u, & x>n
\end{array}\right.
$$
and $b_n(x,u) = \tilde{b}_n(x,u)\ast\zeta_n(x)$ (convolution with
respect to $x$). Similarly, define
$$
\tilde{\sigma}_n(x,u) = \left\{
\begin{array}{ll}
\ds \sigma_0 + \sigma_1 1/n + \sigma_2 u, & |x| < 1/n\\
\ds \sigma_0 + \sigma_1|x| + \sigma_2 u, & 1/n \leq |x| \leq n\\
\ds \sigma_0 + \sigma_1n + \sigma_2 u, & |x|>n,\\
\end{array}\right.
$$
$$
\tilde{\varphi}_n(x) = \left\{
\begin{array}{ll}
\ds \varphi(-n), & x<-n\\
\ds \varphi(x), & -n \leq x \leq n\\
\ds \varphi(n), & x>n,\\
\end{array}\right.
$$
and $\sigma_n(x,u) = \tilde{\sigma}_n(x,u)\ast\zeta_n(x)$,
$\varphi_n(x)=\tilde{\varphi}_n(x)\ast\zeta_n(x)$. Similarly, let
$h_n(x)=h(x)\ast\zeta_n(x)$.
\par\noindent
Given a sequence of functions $f_n(x,u)$, we shall say that $f_n(x,u)$
converges to $f(x,u)$ in $\mathscr{L}$ if for each $R>0$ one has
$$
\lim_{n\to\infty} \sup_{u\in U} \sup_{|x|\leq R}
|f_n(x,u)-f(x,u)| = 0.
$$
By well known properties of convolution with smooth kernels, it is
easy to prove that $\sigma_n$, $b_n$, $h_n$, $\varphi_n$ belong to
$C^3_b$ and that they converge to $\sigma$, $b$, $h$, $\varphi$,
respectively, in $\mathscr{L}$.
Let us denote by $x_t^{u,s,x}(n)$ a solution of the equation
$$
x_t = x + \int_s^t b_n(u_r,x_r)\,dr
+ \int_s^t \sigma_n(u_r,x_r)\,dw_r,
$$
where $u \in \mathfrak{U}$, $0\leq s\leq t \leq T$, $x\in\erre_+$.
Moreover, let us define
\begin{equation}
\label{eq:vn}
v_n^u(s,x) := \E^u_{s,x}\left[
-\int_s^T e^{-cr} h_n(u_r)\,dr + e^{-cT}\varphi_n(x_T(n))
\right]
\end{equation}
and $v_n(s,x)=\sup_{u\in\mathfrak{U}} v_n^u(s,x)$.
\smallskip\par\noindent
\textsc{Step 2:}
Since $\tilde{\sigma}_n(x,u)$ is continuous and
$\tilde\sigma_n^2(x,u)>0$ for all $x\in\erre$ and $u\in U$, by
properties of convolutions with smooth kernels, one also has
$\sigma_n^2(x,u)>0$ uniformly over $x$, $u$ for $n$ large enough.
Therefore (for a fixed large $n$) the value function $v_n$ is a
$C^{1,2}([0,T]\times\erre)$ solution of the Bellman equation
\begin{equation}
\label{eq:ab}
\left\{\begin{array}{l}
\ds \sup_{u\in U}
\Big[ v_t(t,x) + \frac12\sigma_n^2(x,u)v_{xx}(t,x) + b_n(x,u)v_x(t,x)
      - c v(t,x) - h_n(u) \Big] = 0\\[10pt]
v(T,x) = \varphi_n(x),
\end{array}\right.
\end{equation}
as it follows from a result of Krylov (see \cite{krylov-NPDE}, p.~301
and also \cite{FS}, p.~168).
As it is well known, if the supremum in (\ref{eq:ab}) is attained for
each $(x,t)$ by $u^n(x,t)$, then $u^n_t:=u^n(x_t,t)$ is an optimal
Markov strategy for the approximating problem of maximizing
(\ref{eq:vn}). See also \cite{FS}, p.~169.
\smallskip\par\noindent
\textsc{Step 3:}
As it follows by theorem 3.1.12 in \cite{krylov}, we have that
$v_n^u(t,x) \to v^u(t,x)$ uniformly with respect to
$u\in\mathfrak{U}$, $s\in[0,T]$, $|x|\leq R$ for each $R>0$.  Therefore, given
$\varepsilon>0$, there exists $N>0$ such that for all $n>N$ one has
$|v_n^u-v^u|<\varepsilon/3$ for all $u\in\mathfrak{U}$. Let us first
prove that $|v_n-v|<\epsilon$. Assume, by contradiction, that
$v_n<v-\varepsilon$ (without loss of generality, as the case
$v<v_n-\varepsilon$ is completely analogous). By definition of $v$,
there exists $u^1$ such that $|v^{u^1}-v|<\varepsilon/3$. From
$|v_n^{u^1}-v^{u^1}|<\varepsilon/3$ we also have
$$
|v_n^{u^1}-v| \leq |v_n^{u^1}-v^{u^1}|+|v^{u^1}-v| < \frac23\varepsilon,
$$
hence $v_n^{u^1}>v_n$, which is absurd. Then we have proved that
$|v_n-v|<\varepsilon$.
Let us now take $n>N$ large enough, so that the conditions of step 2
are satisfied, and let $u^\epsilon$ be a Markov strategy such that
$v_n=v_n^{u^\epsilon}$. Then we have, by the triangular inequality,
\begin{eqnarray*}
|v^{u^\varepsilon}-v| &\leq& |v^{u^\varepsilon}-v_n^{u^\varepsilon}|
                             + |v_n^{u^\varepsilon}-v| \\
&<& \frac{\varepsilon}{3} + \varepsilon = \frac43\varepsilon,
\end{eqnarray*}
which proves the claim. The proof of the theorem is thus finished.
\end{demo}

Under more specific structural assumptions on the objective function
to optimize, it is natural to expect sharper characterizations of the
value function and of the optimal advertising strategy. This is the
topic of the following sections.


\section{Linear reward and loss functions}
In this section we study the simplest case, with linear reward for the
level of goodwill at time $T$ and linear loss function of advertising
effort. As discussed before, if we identify $u_t$ with the rate of
advertising spending, then the second term in the objective function
(\ref{eq:fo}) can be identified with discounted cumulative advertising
costs.
Let us specify the problem in detail. Our aim is to maximize over
$\mathfrak{U}$ the functional
$$
v^u(s,x) = \E^u_{s,x}\left[
\gamma_0e^{-cT} x_T - \int_0^T e^{-ct} u_t\,dt,
\right]
$$
where
\begin{equation}
\label{eq:sl}
x_t = x + \int_s^t (-\rho x_r + u_r)\,dr + 
\int_s^t \sigma(x_r,u_r)\,dw_r,
\quad s \leq r \leq t \leq T.
\end{equation}
In the sequel we shall set, for simplicity, $\gamma=\gamma_0e^{-cT}$.
Note that, due to the linearity on $u_t$ of the performance functional,
we can explicitly solve the optimization problem even without
specifying the functional form of the diffusion coefficient $\sigma$,
as long as (\ref{eq:sl}) admits a solution.
In particular, $\sigma$ could be identically zero, for which we obtain
the classical deterministic Nerlove-Arrow dynamics.
The interpretation of this fact is simply that optimizing a linear
objective function as $v^u$ simply coincides with controlling the
``mean evolution'' of our stochastic dynamics, for which we can obtain
an explicit expression as follows: first write
$$
x_T = e^{-\rho T}x + \int_0^T e^{-\rho(T-t)}u_t\,dt +
\int_0^T e^{-\rho(T-t)} \sigma(x_t,u_t)\, dw_t,
$$
then take expectations on both sides to get
$$
\E[x_T] = e^{-\rho T}x + \int_0^T e^{-\rho(T-t)}\E[u_t]\,dt,
$$
where the interchange of the order of integration follows by
Fubini's theorem using the assumption $u\geq 0$. It is now easy to
guess that the functional form of $\sigma$ will not influence the
optimal advertising strategy, which is expected to be of the bang-bang
type. This is made precise in what follows. From the managerial point
of view, this means that the firm, independently of the intensity of
the noise and of its dependence on the level of goodwill and rate of
advertising spending, will do its best in terms of maximizing
\emph{expected} goodwill by simply concentrating all its advertising
efforts in a specific period of the advertising campaign.

Assuming $U=[0,m]$, the Bellman equation associated to the problem of
maximizing $v^u$ over $\mathfrak{U}$ is given by
\begin{equation}
\label{eq:hjb-lin}
\psi_t + \sup_{u\in[0,m]}(L^u\psi - e^{-c t}u) = 0,%
\quad \psi(T,x)=\gamma x.
\end{equation}
Note that one has
$$
\sup_{u\in[0,m]}(L^u\psi - e^{-c t}u) = \left\{%
\begin{array}{ll}
-\rho x \psi_x + {1\over 2}\sigma^2 \psi_{xx}, & \psi_x \leq e^{-c t}\cr
-\rho x \psi_x + m(\psi_x-e^{-c t}) + {1\over 2}\sigma^2 \psi_{xx},%
 & \psi_x > e^{-c t}.\cr
\end{array}\right.
$$
Let us consider first the case $\psi_x > e^{-c t}$. Then
(\ref{eq:hjb-lin}) can be written as
$$
\psi_t - (\rho x - m)\psi_x + {1\over 2}\sigma^2\psi_{xx} -me^{-c t} = 0,%
\quad \psi(T,x)=\gamma x.
$$
We guess a solution of the form $\psi(t,x) = \gamma(t)x + b_1(t)$, obtaining
$$ x\gamma'(t) + b_1'(t) - (\rho x - m)\gamma(t) -me^{-c t} = 0, $$
with terminal conditions $\gamma(T)=\gamma$, $b_1(T)=0$.
Then this equation splits into
$$ \gamma'(t)-\rho\gamma(t)=0, \quad \gamma(T)=\gamma, $$
with solution $\gamma(t) = \gamma e^{-\rho(T-t)}$,
and
$$ b_1'(t) = -m\gamma(t) + m e^{-c t}, \quad b_1(T)=0, $$
with solution
$$
b_1(t) = -{m\gamma\over\rho}(1-e^{-\rho(T-t)})+%
{m\over c}(e^{-c T}-e^{-c t}).
$$
The case $\psi_x \leq e^{-c t}$ is completely similar: let $t_*$ be
the solution of the equation $\gamma(t) = e^{-c t}$, i.e.  $t_* =
{\rho T - \log \gamma \over \rho + c}$.  Let us now solve the equation
$$
\psi_t -\rho x \psi_x + {1\over 2}\sigma^2 \psi_{xx} = 0,%
\quad \psi(t_*,x) = \gamma(t_*)x+b_1(t_*),
$$
where the terminal condition is such that a global solution of
(\ref{eq:hjb-lin}) equation is at least continuous.  It is immediate
that the solution of this equation is $\psi(t,x) =
\gamma(t)x+b_1(t_*)$, so that the global solution of
(\ref{eq:hjb-lin}) is $\psi(t,x) = \gamma(t)x+b(t)$, where
$b(t)=b_1(t_*)$ for $\gamma(t) \leq e^{-c t}$, and $b(t)=b_1(t)$ for
$\gamma(t) > e^{-c t}$.  It is also easy to see that $b$ is
continuously differentiable on $(0,T)$.  In fact, one only needs to
check whether there is smooth fit at $t_*$.  But since
$b'_1(t)=-m\gamma(t)+me^{-c t}$, by definition of $t_*$ it immediately
follows $b'_1(t_*)=0$. This also proves that $\psi \in
C^{1,2}([0,T],\erre)$, hence the solution of the Bellman equation
(\ref{eq:hjb-lin}) is the value function of the corresponding control
problem, and we can conclude that the optimal control is given by the
following bang-bang policy:
\begin{equation}
u_*(t) = \left\{\begin{array}{ll}
0 & t \leq t_*, \cr
m & t > t_*. \cr
\end{array}\right.
\label{eq:bang}
\end{equation}
That is, it is optimal not to advertise until a certain point in time $t_*$,
after which it becomes optimal to advertise at the maximum rate. Note that,
depending on $\gamma$, it could well be that $t_*>T$, i.e. it would never
be optimal to advertise. This situation arises if the reward for improving the
image of a product is small compared to the value of resources spent on
advertisement.

We collect the findings of this section in the following proposition.
\begin{prop}
The optimal control problem of maximizing $v^u(0,x)$ is solved by a
control of the type (\ref{eq:bang}), with
$t_*={\rho T-\log\gamma\over\rho+c}$, and the corresponding value
function is given by $v(t,x) = \gamma(t)x+b(t)$, where
$$
b(t) = \left\{\begin{array}{ll}
b_1(t_*) & t \leq t_*, \cr
b_1(t)   & t > t_*. \cr
\end{array}\right.
$$
\end{prop}

Using a Lagrange multiplier method, we can treat the related problem
of maximizing the level of goodwill at time $T$ with a certain
available budget for advertising. More precisely, let us consider the
constrained stochastic control problem
\begin{equation}
\sup_{u\in\mathcal{M}} \E[x_T],
\label{eq:lincon}
\end{equation}
where $\mathcal{M}\subset\mathfrak{U}$ is the set of admissible controls
$u(\cdot)\in[0,m]$ satisfying the integral constraint
$$ 
\E\left[\int_0^T e^{-c t} u_t\,dt\right] \leq M, 
$$
where $M$ is a fixed positive constant. In order for the
constrained problem to be non-trivial, it is also necessary to assume
that $M \leq m\int_0^T e^{-c t}\,dt$. We actually only need to
consider controls $u$ for which the constraint is binding, i.e.
advertising policies that use the whole budget $M$. In fact, denoting
by $x^u_T$ the controlled goodwill at time $T$, it is clear that
$u_1\geq u_2$ implies $\E[x^{u_1}_T] \geq \E[x^{u_2}_T]$, so it is
never optimal to leave resources unused.

Let us introduce a Lagrange multiplier $\lambda>0$, and consider the
(unconstrained) problem
\begin{equation}
\sup_{u\in[0,m]} \E\left[x_T -%
\lambda\left(\int_0^T e^{-c t} u_t\,dt-M\right)\right].
\label{unc-prob}
\end{equation}
Then one has
\begin{eqnarray*}
\sup_{u\in\mathcal{M}} \E[x_T] &=&
\sup_{u\in\mathcal{M}} \E\left[x_T -%
\lambda\left(\int_0^T e^{-c t} u_t\,dt-M\right)\right] \\
&\leq& \sup_{u\in\mathfrak{U}} \E\left[x_T -%
\lambda\left(\int_0^T e^{-c t} u_t\,dt-M\right)\right],
\end{eqnarray*}
where the first equality comes from the above observation that we only
need consider controls that use the whole budget $M$, and the second
inequality is justified by $\mathcal{M} \subseteq \mathfrak{U}$.
If the unconstrained problem (\ref{unc-prob}) admits a solution $u_\lambda$
for all $\lambda>0$, and a $\lambda_*$ exists such that
$\E\int_0^T e^{-c t} u_{\lambda_*}(t)\,dt-M=0$,
then $u_*:=u_{\lambda_*}$ is an optimal control for the constrained problem.
So we proceed to solve
$$
\sup_{u\in\mathfrak{U}} \E\left[{1\over\lambda} x_T - \int_0^T e^{-c t}
u_t\,dt\right],
$$
whose solution is
\begin{eqnarray*}
\gamma(t) \leq e^{-c t} &\Rightarrow& u_*(t) = 0, \\
\gamma(t) >    e^{-c t} &\Rightarrow& u_*(t) = m,
\end{eqnarray*}
with $\ds \gamma(t) = {1\over\lambda}e^{-\rho(T-t)}.$

The starting point for advertisement $t_*$ is given by the solution of the
equation $\gamma(t) = e^{-c t}$,
so that
\begin{equation}
t_* = {\rho T + \log\lambda \over \rho + c}.
\label{eq:tlambda}
\end{equation}
We now need to show that $\lambda > 0$ exists such that
\begin{equation}
\int_{t_*}^T m e^{-c t}\,dt = M.
\label{eq:lambda}
\end{equation}
The solution of such an equation is given by
$$
\lambda_* = e^{\rho T} \Big(c{M\over m}+%
e^{-c T}\Big)^{-{\rho+c\over c}},
$$
which is clearly positive. It is now clear how to associate to such a
$\lambda_*$ the optimal solution for the constrained problem. Namely, given
$\lambda_*$ we obtain the optimal switching time $t_*$ by
(\ref{eq:tlambda}), and hence the optimal control as
$u_*(t) = m\chi_{\{t>t_*\}}$, where $\chi$ is the indicator
function.

We have then proved the following result.
\begin{prop} The optimal advertising policy for the constrained maximization
of goodwill (\ref{eq:lincon}) is given by
$$
u_*(t) = \left\{\begin{array}{ll}
0 & t \leq t_*, \cr
m & t > t_*, \cr
\end{array}\right.
$$
with
$$ t_* = {2\rho\over\rho+c}T - {1\over c}(e^{-c T} + c M/m).$$
\end{prop}
\begin{rmk}
It follows from (\ref{eq:lambda}) that the time to start
advertising is given by
$$ e^{-c t_*} - e^{-c T} = c{M\over m}, $$
and therefore we cannot simply consider unbounded controls with cumulative
discounted cost less or equal than $M$, otherwise the optimal policy would
be to ``do infinite advertising at time $T$''.
\end{rmk}


\section{Quadratic reward and loss functions}
In this section we assume that both $\varphi$ and $h$ are quadratic. In
particular, we assume $\varphi(x)=\gamma_0 x^2$, with $\gamma_0>0$, and
$h(x)=x^2$. That is, we consider the problem of characterizing
$$ 
v(s,x)= \sup_{u\in\mathfrak{U}} \E^u_{s,x}\Big[ \gamma x_T^2 - \int_0^T e^{-c
t} u_t^2 \,dt \Big],
$$
or, equivalently,
\begin{equation}
v(s,x) = \inf_{u\in\mathfrak{U}} \E^u_{s,x}\Big[ \int_0^T
e^{-c t} u_t^2 \,dt - \gamma x_T^2 \Big],
\label{LQ-prob}
\end{equation}
where we set $\gamma=e^{-cT}\gamma_0$, $x_t$ follows the controlled
dynamics
$$
x_t = x + \int_s^t (-\rho x_r+u_r)\,dr 
        + \int_s^t (\sigma_1 x_r + \sigma_2 u_r)\,dw_r,
$$
and $\mathfrak{U}$ is the set of adapted, nonnegative, square
integrable controls.
A peculiar feature of the problem is that, while the choice of the
cost function $h$ is rather standard (see e.g. Muller \cite{muller}),
the reward function $\varphi$ is convex, hence representative of a
risk-seeking firm. Such attitude toward risk could be justified, for
instance, by the attempt to profit from the fluctuations of the
goodwill level at time $T$ (in fact, note that, grossly speaking, the
firm aims at maximizing both the mean and the variance of $x_T$).
Another peculiar feature of $\varphi(x)=x^2$ is that it equally
rewards positive and negative goodwill levels at time $T$. However, we
shall show that under our assumptions $x_t \geq 0$ almost surely,
hence the symmetry of $\varphi$ is harmless.

In this section we assume that the intensity of the ``background
noise'' is negligible, so that we can assume $\sigma_0=0$. This
assumption is essential in order to obtain (meaningful) solutions in
closed form.

The problem at hand is a linear quadratic regulator problem with
indefinite costs, which can be solved by the methods of Ait Rami,
Moore and Zhou \cite{AMZ} (see also Krylov \cite{krylov-LQ}). In
particular, the generalized Riccati equation for this problem is
\begin{equation}
\label{eq:R3}
\left\{\begin{array}{l}
\ds\dot{P} = (2\rho-\sigma_1^2)P + %
(1+\sigma_1\sigma_2)^2 {P^2 \over e^{-c t}+\sigma_2^2P}, \cr
e^{-c t}+\sigma_2^2P > 0 \cr
P(T)=-\gamma.
\end{array}\right.
\end{equation}
Recall that one says that (\ref{LQ-prob}) is well posed at $s$ if
$v(s,x)>-\infty$. By \cite{AMZ}, well-posedeness of (\ref{LQ-prob}) at
$s=0$ is necessary for the global solvability of (\ref{eq:R3}).
Under the assumption $e^{-c T}-\sigma_2^2\gamma>0$, one can only
ensure that the problem is locally well posed, i.e. that there exists
$t_0<T$ such that the Riccati equation (\ref{eq:R3}) admits a solution
in $[t_0,T]$. The following proposition gives explicit sufficient
conditions on the data of the problem such that (\ref{eq:R3}) admits a
unique global solution on $[0,T]$.
\begin{prop}
Let us define
\begin{eqnarray*}
a_1 &=& -2\rho - 2\sigma_1\sigma_2^{-1} - \sigma_2^{-2} < 0\\
a_2 &=& 2\rho + \sigma_1^2 + 2\sigma_2^{-2} + 4\sigma_1\sigma_2^{-1} > 0\\
a_3 &=& -(\sigma_1 + \sigma_2^{-1})^2 < 0\\
a_4 &=& -\gamma_0\sigma_2^2 + 1 > 0\\
\zeta &=& a_2^2-4a_1a_3 = 4\rho^2 + \sigma_1^4 - 4\rho\sigma_1^2,
\end{eqnarray*}
where the inequality $a_4>0$ is taken as an assumption. Then the
following hold:
\begin{itemize}
\item[(i)] If $\zeta>0$ and $a_2 > (2|a_1|a_4 - \zeta^{1/2})^+$, then
  the problem is well posed.
\item[(ii)] If $\zeta>0$ and $a_2 \leq (2|a_1|a_4 - \zeta^{1/2})^+$,
  then the problem is well posed if and only if
  $$
  T \leq \zeta^{-1/2}\Big( \xi_1\log\frac{\xi_1}{\xi_1-a_4} -
  \xi_2\log\frac{\xi_2}{\xi_2-a_4} \Big),
  $$
  where $\xi_{1,2} = (-a_2 \pm \zeta^{1/2})/(2a_1)$.
\item[(iii)] If $\zeta=0$ and $a_2>2|a_1|a_4$, then the problem is well
  posed.
\item[(iv)] If $\zeta=0$ and $a_2\leq 2|a_1|a_4$, then the problem is well
  posed if and only if
  $$
  T \leq a_1^{-1}\Big( \log\frac{a_2}{2a_1a_4+a_2} +
  \frac{2a_1a_4}{2a_1a_4+a_2} \Big)
  $$
\item[(v)] If $\zeta<0$, then the problem is well posed if and only if
$$
\begin{array}{l}
T \leq (2a_1)^{-1}\bigg(
\log\frac{a_3}{a_1a_4^2+a_2a_4+a_3} \\
\phantom{T \leq (2a_1)^{-1}\bigg(}
-2a_2\zeta^{-1/2}\Big(\mathrm{atan}\,a_2\zeta^{-1/2}
-\mathrm{atan}\,\zeta^{-1/2}(2a_1a_4+a_2)\Big)
\bigg)
\end{array}
$$
\end{itemize}
\end{prop}
\begin{demo}
  We shall divide the proof in three steps. In the first step we
  introduce an auxiliary LQ problem whose well-posedness is sufficient
  for the well-posedness of (\ref{LQ-prob}). In the second step we
  rescale the auxiliary LQ problem, and in the third and last step we
  study the global solvability of its associated Riccati equation.
\par\noindent
\textsc{Step 1:} Let us prove that if
\begin{equation}
\label{eq:LQ-aux}
\inf_{u\in\mathfrak{U}} \E\left[
\int_0^T u_t^2\,dt - \gamma_0 x_T^2 \right] > -\infty,
\end{equation}
then (\ref{LQ-prob}) is well posed. In fact, suppose, by
contradiction, that there exists a sequence $u(k) \in \mathfrak{U}$
such that
$$
\E^{u(k)}\left[
\int_0^T e^{-ct} u_t(k)^2\,dt - \gamma x_T^2 \right] \to -\infty.
$$
Then one also has
\begin{eqnarray*}
e^{-cT} \E^{u(k)}\left[
\int_0^T u_t(k)^2\,dt - \gamma_0 x_T^2 \right] &\leq&
\E^{u(k)}\left[
\int_0^T e^{-ct}u_t(k)^2\,dt - \gamma x_T^2 \right] \to -\infty,
\end{eqnarray*}
which contradicts (\ref{eq:LQ-aux}), hence our claim is proved.
\smallskip\par\noindent
\textsc{Step 2:} Let us define $\tilde{x}_t=\sigma_2^{-1}x_t$. Then
the minimization problem in (\ref{eq:LQ-aux}) is equivalent to
\begin{equation}
\label{eq:LQa}
\inf_{u\in\mathfrak{U}} \E\left[
\int_0^T u_t^2\,dt - \gamma_0\sigma_2^2 \tilde{x}_T^2 \right]
\end{equation}
subject to
$$
\tilde{x}_t = \sigma_2^{-1}x 
+ \int_0^t (-\rho\tilde{x}_s+\sigma_2^{-1}u_s)\,ds 
+ \int_0^t (\sigma_1\tilde{x}_s+u_s)\,dw_s.
$$
\smallskip\par\noindent 
\textsc{Step 3:} The Riccati equation for problem (\ref{eq:LQa}) is
$$
\left\{
\begin{array}{l}
\ds \dot{P} = (2\rho-\sigma_1^2)P +
\frac{(\sigma_1+\sigma_2^{-1})^2P^2}{P+1} \\[6pt]
P(T) = -\gamma_0\sigma_2^2 \\[6pt]
P(t)+1 > 0,
\end{array}\right.
$$
which can be rewritten, after the change of variable
$\pi(t)=P(T-t)+1$, as
\begin{equation}
\label{eq:rr}
\left\{
\begin{array}{l}
\ds \dot{\pi} = a_1\pi + a_2 + a_3\pi^{-1} \\
\ds \pi(0) = a_4 \\
\ds \pi(t) > 0.
\end{array}\right.
\end{equation}
All assertions (i)-(v) are proved simply by determining the first time
$t_0$ for which $\pi(t_0)=0$. Since the method of proof is the same
for all cases, we shall show in detail only cases (i)-(ii).  In
particular, assuming $\zeta>0$, then $a_1\pi^2+a_2\pi+a_3=0$ has two
roots $\xi_1<\xi_2$, and the reduced Riccati equation (\ref{eq:rr}) as
$$
\dot{\pi} = a_1\pi^{-1} (\pi-\xi_1)(\pi-\xi_2).
$$
Note that $a_1<0$ implies $\xi_1\xi_2=a_3/a_1>0$. Moreover one has
$\xi_1+\xi_2=a_2/|a_1|$, and, by $a_2>0$, $0<\xi_1<\xi_2$. Therefore
$t_0<\infty$ if and only if $a_4<\xi_1 =
\frac{-a_2-\zeta^{1/2}}{2a_1}$, which can be rewritten as
$a_2>(2|a_1|a_4-\zeta^{1/2})^+$.
In order to determine $t_0$ such that $\pi(t_0)=0$, we shall solve
(\ref{eq:rr}) explicitly: for this purpose, note that one has
$$
\frac{\pi}{(\pi-\xi_1)(\pi-\xi_2)} = \frac{A}{\pi-\xi_1} +
\frac{B}{\pi-\xi_2},
$$
with $A=-\frac{\xi_1}{\xi_2-\xi_1}$ and $B=\frac{\xi_2}{\xi_2-\xi_1}$, hence
(\ref{eq:rr}) can also be written as
$$
A\frac{\dot{\pi}}{\pi-\xi_1} + B\frac{\dot{\pi}}{\pi-\xi_2} = a_1.
$$
Integrating one gets
$$
-\xi_1 \int_0^t \frac{\dot{\pi}(s)}{\pi(s)-\xi_1}\,ds +
\xi_2 \int_0^t \frac{\dot{\pi}(s)}{\pi(s)-\xi_2} = a_1(\xi_2-\xi_1)t,
$$
and finally 
$$
-\xi_1 \log\frac{\pi(t)-\xi_1}{\pi(0)-\xi_1} +
\xi_2 \log\frac{\pi(t)-\xi_2}{\pi(0)-\xi_2} = a_1(\xi_2-\xi_1)t.
$$
Solving for $t_0$ such that $\pi(t_0)=0$, recalling that
$a_1(\xi_2-\xi_1)=\zeta^{1/2}$, gives the required result.
\end{demo}

Assuming that the problem is well posed, the results in \cite{AMZ} imply that
the optimal control strategy is unique and is given by the Markov policy
\begin{equation}
\label{eq:lqo}
u^0_t = u^0(x_t) =
-{ (1+\sigma_1\sigma_2)P(t) \over e^{-c t}+\sigma_2^2P(t)} x_t, 
\end{equation}
with associated value function $v(s,x)=P(s)x^2$.
The optimal trajectory is then given by the closed-loop equation
$$ dx_t = a(t)x_t\,dt + c(t)x_t\,dw_t, $$
with
\begin{eqnarray*}
a(t) &:=& -\rho - {(1+\sigma_1\sigma_2)P(t) \over e^{-c t}+\sigma_2^2P(t)}, \\
c(t) &:=& \sigma_1 + \sigma_2{(1+\sigma_1\sigma_2)P(t) \over e^{-c t}+\sigma_2^2P(t)},
\end{eqnarray*}
which admits the explicit solution
\begin{equation}
\label{eq:traj3}
x_t = 
x \exp\Big(\int_0^t(a(s)-{1\over 2}c(s)^2)\,ds + \int_0^t c(s)\,dw_s\Big).
\end{equation}
In particular, if $x>0$, then $x_t$ is a.s. positive for all
$t\in[0,T]$. If we can prove that $P(t)<0$ for all $t\in[0,T]$, then
(\ref{eq:lqo}) will imply that the optimal strategy $u^0_t$ is
positive for all $t\in[0,T]$. The negativity of $P$ is proved in the
following lemma.

\begin{lemma} 
  If $P$ solves the Riccati equation (\ref{eq:R3}) on $[s,T]$, then
  $P(t)<0$ for all $t \in [s,T]$.
\end{lemma}
\begin{demo} 
  Assume, by contradiction, that there exists $t_0\in[s,T[$ such that
  $P(t_0)=0$. Then $P(t) \equiv 0$ is a continuous solution of
  (\ref{eq:R3}) for $t>t_0$. By uniqueness of the solution of
  (\ref{eq:R3}) (see \cite{AMZ}) it follows that this is also the only
  solution. But this contradicts the terminal condition $P(T)=-\gamma
  \neq 0$.
\end{demo}
We have thus proved all claims on the optimal state and control.
\smallskip\par\noindent
\textbf{Remarks.}
  (i) It is worth noting that several sensitivities of the value function
  and of the expected optimal goodwill with respect to initial data or
  parameters could be computed as well in terms of the solution of the
  Riccati equation (\ref{eq:R3}). 

  (ii) If the intensity of the noise carried by the advertising is also
  negligible, i.e. we can assume $\sigma_2=0$, then the Riccati
  equation (\ref{eq:R3}) is explicitly solvable, and the corresponding
  calculations mentioned in the previous remark simplify even more.
  More details can be found in \cite{gw}.
  
  (iii) The linear-quadratic regulator approach is also useful if one
  wants to consider problems with partial observation, which are
  meaningful in our setting, as goodwill can hardly be measured
  without error. Suppose instead that we can observe a ``noisy
  proxy'' of goodwill $z_t$ specified by
\begin{equation}
\label{eq:observ}
\left\{\begin{array}{rcl}
dz_t \!&=&\! h x_t\,dt + g\,dw^1_t \cr
z_0 \!&=&\! 0,
\end{array}\right.
\end{equation}
with $h$ and $g$ constants, and $w^1$ a Brownian motion independent of
$w$.
Then one is interested to solve problem (\ref{LQ-prob}), where
$\mathfrak{U}$ is now the set of nonnegative square integrable
controls adapted to the filtration generated by $z_t$, instead of $x_t$.
Thanks to the so called separation principle (see e.g. Bensoussan
\cite{B}), this problem reduces to linear filtering and deterministic
control on the filtered dynamics. More details will be given
elsewhere.
\smallskip


\section{Optimal advertising to meet a goal with discretionary stopping}
Problems of mixed optimal stopping and control have recently attracted
attention in works of applied probability, see for instance Karatzas and
Wang \cite{KW} for applications to portfolio optimization, Duckworth and
Zervos \cite{duck1}, \cite{duck2} and Zervos \cite{zervos} for problems of
investment decisions with strategic entry and exit, and Karatzas, Ocone,
Wang and Zervos \cite{KOWZ} for a singular control problem with finite fuel. 
For the theory, see, e.g., Krylov \cite{krylov}, Bensoussan and Lions
\cite{BL}, {\O}ksendal and Sulem \cite{OS}, and Morimoto \cite{morimoto}. 
One of the first works addressing the issue of finding explicit results was
Bene\v{s} \cite{benes}.

In this section we find an explicit representation for the optimal control
and the optimal stopping strategy for the case of minimizing an objective
function that is the sum of the quadratic distance of the goodwill from a
target at the (discretionary) launch time $\tau$ and of the cumulative
quadratic cost until $\tau$, assuming that the goodwill dynamics is of
Ornstein-Uhlenbeck type.
For simplicity we also assume $c=0$, i.e. we consider the case without
discounting. In order to discourage long waiting before launching the
product, we also introduce an extra term in the objective function depending
on the time of launching.

Let $y_t=k-x_t$ be the distance of the goodwill level at time $t$ from
a desired target $k$. We shall find the solution to the problem
\begin{equation}
\inf_{u\in\mathfrak{U},\tau\in\mathcal{S}}\;
\E^{u,\tau}_{0,x}\!\left[
y^2(\tau) + \gamma_1\int_0^\tau u^2(t)\,dt + \gamma_2\tau
\right]=:v(x),
\label{eq:stopc}
\end{equation}
where $y$ is such that
$$ dy_t=(\mu-\rho y_t+u_t)\,dt + dw_t, $$
with $\mu:=\rho k$ and we have assumed, without loss of generality (in the
setting of constant $\sigma$), 
$\sigma=1$. Let $\mathbb{F}$ be the filtration generated by $w$. Then
$\mathfrak{U}$ is the space of $\mathbb{F}$-adapted square integrable control
processes, and $\mathcal{S}$ is the set of all $\mathbb{F}$-stopping times.

Note that the first term in (\ref{eq:stopc}) assigns equal costs to
the events that the target is not reached (from below) and that it is
exceeded. Our setting could be considered as a stylized model for the
situation when a firm is launching a new product in a predetermined
quantity (the target $k$) and considers equally undesirable to
undersell the product ($x_t<k$) or to leave demand unmet ($x_t>k$).

The quasi-variational inequality associated to the mixed problem of optimal
control and optimal stopping (\ref{eq:stopc}) can be written as
$$ 
\min_x\left(x^2-v(x),\min_u(L^uv+\gamma_1 u^2 + \gamma_2)\right) = 0,
$$
where $L^u$ is the generator of the controlled diffusion $y$, i.e.
$L^u$ is the differential operator defined by
$$ L^u f(x) = {1\over 2}f''(x)+(\mu-\rho x-u)f'(x). $$
We guess a continuation region $D$ of the type $D=\{x:x\geq x_0\}$, where one
must have
$$ \min_u(L^uv+\gamma_1 u^2 + \gamma_2) = 0. $$
We have
$$
\min_u(L^uv+\gamma_1 u^2 + \gamma_2) = Av - {1\over 4\gamma_1}v_x^2 +
\gamma_2,
$$
where $A$ is the generator of the uncontrolled diffusion, i.e.
$$ A f(x) = {1\over 2}f''(x)+(\mu-\rho x)f'(x). $$
Then we get
$$ Av - {1\over 4\gamma_1}v_x^2 + \gamma_2 = 0, \quad x\geq x_0 $$
In order to linearize this ODE, we apply the Hopf-Cole
transformation $U(x)=e^{{1\over 2\gamma_1}v(x)}$, obtaining
$$
{1\over 2}U_{xx} +(\mu-\rho x)U_x - {\gamma_2\over 2 \gamma_1}U=0, %
\quad x\geq x_0
$$
In order to obtain solutions that are ordinary functions, we restrict
ourselves to the special case ${\gamma_2\over 2\gamma_1}=\rho$.
However, one can solve the linear equation for $U$ without this
assumption, obtaining solutions that can be expressed in terms of
special functions.
Two linearly independent solutions are
$$
U_1(x) = e^{-2\mu x+\rho x^2}, \quad
U_2(x) = e^{\rho(x-\mu/\rho)^2}\int_x^\infty e^{-\rho(s-\mu/\rho)^2}\,ds
$$
Then the general solution can be written as $U=\alpha_1 U_1 +
\alpha_2 U_2$, with $\alpha_1$ and $\alpha_2$ arbitrary real
constants.

Guided by the observation that $U_1$ is unbounded in the continuation
region, we set $\alpha_1=0$, and impose $C^1$ fit of
$U(x)=\alpha_2U_2(x)$ to the Hopf-Cole transformation of $x^2$ at the
point $x_0$, that is
\begin{eqnarray*}
U(x_0) &=& e^{-{1\over2\gamma_1}x_0^2} \\
U'(x_0) &=& -{x_0\over\gamma_1} e^{-{1\over2\gamma_1}x_0^2}.
\end{eqnarray*}
We solve now the following system of equations for the unknowns
$\alpha_2$ and $x_0$:
\begin{eqnarray*}
\alpha_2 e^{\rho(x_0-\mu/\rho)^2}\int_x^\infty e^{-\rho(s-\mu/\rho)^2}\,ds &=& e^{-{1\over2\gamma_1}x_0^2} \\
\alpha_2 \Big[2\rho(x_0-{\mu\over\rho})e^{\rho(x_0-\mu/\rho)^2}%
\int_{x_0}^\infty e^{-\rho(s-\mu/\rho)^2}\,ds - 1\Big] &=&
-{x_0\over\gamma_1} e^{-{1\over2\gamma_1}x_0^2}.
\end{eqnarray*}
Therefore $x_0$ is given by the solution of the following equation
$$
((2\rho+{1\over\gamma_1})x-2\mu)%
e^{\rho(x-\mu/\rho)^2}%
\int_x^\infty e^{-\rho(s-\mu/\rho)^2}\,ds = 1.
$$
In fact the solution, if it exists, is unique, because the left hand side, as a function of $x$, is increasing.
Now one can also find $\alpha_2$ in terms of $x_0$, so that we have a
candidate continuation region (or equivalently an optimal stopping
region), and a candidate value function.

We need to show that $Av(x)-{1\over 4\gamma_1}(v'(x))^2+\gamma_2\geq0$ in the region $x\leq x_0$, for $v(x)=x^2$.
That is, we want to show that
\begin{equation}
(2\rho+{1\over\gamma_1})x^2 - 2\mu x - (1+\gamma_2) \leq 0.
\label{eq:cond1}
\end{equation}
The expression on the left hand side takes its maximum either at $x=0$ or at
$x=x_0$. Therefore, if $x_{max}=0$, condition (\ref{eq:cond1}) is trivially
verified. Otherwise, if $x_{max}=x_0$, we have
$$
{1 \over (2\rho+\gamma_1^{-1})x_0-2\mu} = U_2(x_0) \geq
{1 \over \rho y_0 + \sqrt{\rho^2y_0^2+2\rho}},
$$
where $y_0:=x_0-{\mu\over\rho}$, and the inequality follows by standard
estimates on the error function.
Therefore we get
$$
(2\rho+\gamma_1^{-1})x_0-2\mu \leq
\rho y_0 + \sqrt{\rho^2y_0^2+2\rho}.
$$
After some algebraic manipulations, one finds
$$
(2\rho+{1\over\gamma_1})x_0^2 - 2\mu{1\over\gamma_1}x_0
\leq 2\rho\gamma_1 = \gamma_2.
$$
Therefore, (\ref{eq:cond1}) is verified if, e.g., $\gamma_1>1$. We
shall assume that in the following, but note that this condition can
be weakened.

Finally, we need to prove that $v(x) \geq x^2$ in the continuation region $x\geq x_0$, or equivalently that
$U(x) \geq e^{-{1\over2\gamma_1}x^2}$ for $x\geq x_0$. In order to prove
this, let us consider their ratio
$f(x) = U(x)e^{{1\over2\gamma_1}x^2}$ and prove that it is increasing.
One has
$$
f'(x) = \alpha_2e^{{1\over2\gamma_1}x^2}%
\Big( ((2\rho+{1\over\gamma_1})x-2\mu)U_2(x)-1 \Big),
$$
and since we have that $((2\rho+{1\over\gamma_1})x-2\mu)U_2(x)$ is
increasing in $x$ and
$((2\rho+{1\over\gamma_1})x_0-2\mu)U_2(x)=1$, it follows
$((2\rho+{1\over\gamma_1})x-2\mu)U_2(x)>1$ for $x\geq x_0$.
Therefore we have verified all conditions for optimality, and we
summarize our findings in the following proposition.

\begin{prop} The optimal control policy $u_*$ and optimal stopping time
$\tau_*$ for the problem (\ref{eq:stopc}), with $\gamma_1>1$ and 
${\gamma_2\over2\gamma_1}=\rho$, are given by
$$ 
u_*(y_t) = \mathrm{arg}\min_u (L^uv(y_t)+\gamma_1 u^2 + \gamma_2) 
          = {v'(y_t) \over 2\gamma_1}
$$
and
$$ \tau_* = \inf\{t\geq 0:\,y_t\leq x_0\}. $$
\end{prop}

\section{Further problems}
The optimal control problems studied in this paper are limited to the case
of ``smooth'' disturbances, that is, the driving noise process has
continuous paths. It is meaningful to relax this assumption and consider
also jump components in the noise, to take into account possible shocks to
the image of the advertised product, due, for instance, to bad news on the
product itself or similar ones, or to the introduction of superior
technologies.

One could also try to study different type of controls, namely impulse
controls, or even combinations of classical and impulse controls. This is
particularly meaningful for our problems, since impulse controls correspond
to the so-called ``pulsing advertising'' policies that have been studied in
the management and marketing literature (see \cite{sethi} and references
therein).

\subsection*{Acknowledgements}
The author wishes to thank Sergio Albeverio, Victor de la Pe\~na,
Fausto Gozzi, Cristian Pasarica, Sergei Savin, and Luciano Tubaro for
helpful comments and suggestions on an earlier version of this paper.
The comments of two anonymous referees, to whom the author is
grateful, considerably improved the presentation of the paper. The
financial support of the National Science Foundation under grant
DMS-02-05791 (Principal Investigator V. de la Pe\~na) is gratefully
acknowledged. Large part of this work was carried out at the Graduate
Business School of Columbia University.

\def\cprime{$'$}


\begin{thebibliography}{10}

\bibitem{AMZ}
M.~Ait~Rami, J.~B. Moore, and X.~Y. Zhou.
\newblock Indefinite stochastic linear quadratic control and generalized
  differential {R}iccati equation.
\newblock {\em SIAM J. Control Optim.}, 40(4):1296--1311, 2001/02.

\bibitem{bass}
F.~M. Bass.
\newblock A new product growth model for consumer durables.
\newblock {\em Manag. Sci.}, 15:215--227, 1969.

\bibitem{benes}
V.~E. Bene\v{s}.
\newblock Some combined control and stopping problems.
\newblock Unpublished manuscript, 1993. Available from the author.

\bibitem{B}
A.~Bensoussan.
\newblock {\em Stochastic control of partially observable systems}.
\newblock Cambridge UP, 1992.

\bibitem{BL}
A.~Bensoussan and J.-L. Lions.
\newblock {\em Contr\^ole impulsionnel et in\'equations quasi variationnelles}.
\newblock Gauthier-Villars, Paris, 1982.

\bibitem{BG}
A.~Buratto and L.~Grosset.
\newblock A communication mix for an event planning: a stochastic approach.
\newblock {\em CEJOR Cent. Eur. J. Oper. Res.}, forthcoming.

\bibitem{BV}
A.~Buratto and B.~Viscolani.
\newblock New product introduction: goodwill, time and advertising cost.
\newblock {\em Math. Methods Oper. Res.}, 55(1):55--68, 2002.

\bibitem{dube}
J.~P. Dube and P.~Manchanda.
\newblock Differences in dynamic brand competition across markets: An empirical
  analysis.
\newblock {\em Marketing Science}, forthcoming.

\bibitem{duck1}
K.~Duckworth and M.~Zervos.
\newblock An investment model with entry and exit decisions.
\newblock {\em J. Appl. Probab.}, 37(2):547--559, 2000.

\bibitem{duck2}
K.~Duckworth and M.~Zervos.
\newblock A model for investment decisions with switching costs.
\newblock {\em Ann. Appl. Probab.}, 11(1):239--260, 2001.

\bibitem{sethi}
G.~Feichtinger, R.~Hartl, and S.~Sethi.
\newblock Dynamical {O}ptimal {C}ontrol {M}odels in {A}dvertising: {R}ecent
  {D}evelopments.
\newblock {\em Management Sci.}, 40:195--226, 1994.

\bibitem{FS}
W.~H. Fleming and H.~M. Soner.
\newblock {\em Controlled {M}arkov processes and viscosity solutions}.
\newblock Springer-Verlag, New York, 1993.

\bibitem{GV}
L.~Grosset and B.~Viscolani.
\newblock Advertising for a new product introduction: a stochastic approach.
\newblock {\em Top}, 12(1):149--167, 2004.

\bibitem{KOWZ}
I.~Karatzas, D.~Ocone, H.~Wang, and M.~Zervos.
\newblock Finite-fuel singular control with discretionary stopping.
\newblock {\em Stochastics Stochastics Rep.}, 71(1-2):1--50, 2000.

\bibitem{KW}
I.~Karatzas and H.~Wang.
\newblock Utility maximization with discretionary stopping.
\newblock {\em SIAM J. Control Optim.}, 39(1):306--329, 2000.

\bibitem{krylov}
N.~V. Krylov.
\newblock {\em Controlled diffusion processes}.
\newblock Springer-Verlag, New York, 1980.

\bibitem{krylov-NPDE}
N.~V. Krylov.
\newblock {\em Nonlinear elliptic and parabolic equations of the second order}.
\newblock ``Nauka'', Moscow, 1985.
\newblock In Russian. English translation: D. Reidel Publishing Co., Dordrecht,
  1987.

\bibitem{krylov-LQ}
N.~V. Krylov.
\newblock Stochastic linear controlled systems with quadratic cost revisited.
\newblock In {\em Stochastics in finite and infinite dimensions}. Birkh\"auser,
  Boston, 2001.

\bibitem{gw}
C.~Marinelli.
\newblock The stochastic goodwill problem.
\newblock \texttt{arXiv:math.OC/0310316} preprint, 2003.

\bibitem{morimoto}
H.~Morimoto.
\newblock Variational inequalities for combined control and stopping.
\newblock {\em SIAM J. Control Optim.}, 42(2):686--708, 2003.

\bibitem{muller}
E.~Muller.
\newblock Trial/awareness advertising decision: A control problem with phase
  diagrams with non-stationary boundaries.
\newblock {\em J. Econ. Dynamics Control}, 6:333--350, 1983.

\bibitem{NA}
M.~Nerlove and J.~K. Arrow.
\newblock Optimal advertising policy under dynamic conditions.
\newblock {\em Economica}, 29:129--142, 1962.

\bibitem{OS}
B.~{\O}ksendal and A.~Sulem.
\newblock {\em Applied stochastic control of jump diffusions}.
\newblock Springer-Verlag, Berlin, 2005.

\bibitem{raman}
K.~Raman.
\newblock Stochastically optimal advertising policies under dynamic conditions:
  the ratio rule.
\newblock {\em Optimal Control Appl. Meth.}, 11:283--288, 1990.

\bibitem{rao}
R.~C. Rao.
\newblock Estimating continuous time advertising-sales models.
\newblock {\em Management Sci.}, 5(2):125--142, 1986.

\bibitem{tapiero}
C.~S. Tapiero.
\newblock {\em Applied Stochastic Models and Control in Management}.
\newblock North Holland, Amsterdam, 1988.

\bibitem{vilca}
N.~J. Vilcassim, V.~Kadiyali, and P.~K. Chintagunta.
\newblock Investigating dynamic multifirm market interactions in price and
  advertising.
\newblock {\em Management Sci.}, 45:499--518, 1999.

\bibitem{zeleny}
M.~Zeleny.
\newblock {\em Multiple Criteria Decision Making}.
\newblock McGraw-Hill, New York, 1981.

\bibitem{zervos}
M.~Zervos.
\newblock A problem of sequential entry and exit decisions combined with
  discretionary stopping.
\newblock {\em SIAM J. Control Optim.}, 42(2):397--421, 2003.

\end{thebibliography}
\end{document}